\documentclass[12pt]{amsart}

\usepackage{a4wide}

\usepackage{epsfig}

\usepackage{boxedminipage}

\usepackage{amssymb}

\newtheorem{thm}[subsection]{Theorem}

\newtheorem{lemma}[subsection]{Lemma}

\newtheorem{coro}[subsection]{Corollary}

\newtheorem{prop}[subsection]{Proposition}

\newtheorem*{maintheorem}{Theorem~\ref{thm:main}}
\newtheorem*{mainlemma}{Lemma~\ref{lemma:bound}}

\theoremstyle{definition}

\newtheorem{remark}[subsection]{Remark}
\newtheorem{question}[subsection]{Question}

\title{The minimum dilatation of pseudo-Anosov 5-braids} 
\author{Ji-Young Ham \and Won Taek Song}

\address{Department of Mathematics,  
University of California Santa Barbara, CA 93106}
\email{ham@mail.math.ucsb.edu}

\address{School of Mathematics,
Korea Institute for Advanced Study,
207-43 Cheongnyangni 2-dong, Dongdaemun-gu,
Seoul 130-722, Korea}
\email{cape@kias.re.kr}


\date{14 June 2005; Revised 24 January 2006}

\begin{document}
\begin{abstract}
The minimum dilatation of pseudo-Anosov 5-braids 
is shown to be the largest zero $\lambda_5 \approx 1.72208$ 
of $x^4 - x^3 - x^2 - x + 1$ which is attained by  
$\sigma_1\sigma_2\sigma_3\sigma_4\sigma_1\sigma_2$.
\end{abstract}
\maketitle

\section{Introduction}
\begin{figure}
$$
\epsfig{file=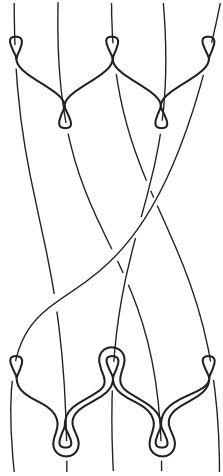,height=0.3\textheight}
$$
\caption{This picture shows a pseudo-Anosov 5-braid 
$\sigma_3 \sigma_4 \sigma_3 \sigma_2 
\sigma_3\sigma_1$ with its invariant train track.  
The train track $\tau \subset D^2 \times 1 $ on the top 
slides down in the complement of the braid 
to $f(\tau) \subset D^2 \times 0$ on the bottom, 
where $f$ is the disk homeomorphism representing  
the braid.}\label{fig:horseshoe}
\end{figure}

Let $f\colon D^2\to D^2$ be an orientation preserving 
disk homeomorphism which 
is the identity map on the boundary $\partial D^2$, 
and $\{p_i\} \subset int(D^2)$ be  a periodic orbit of~$f$ 
(or more generally a finite set invariant under~$f$). 
The points $p_i$ move under an isotopy 
from the identity map on $D^2$ to $f$. 
Their trajectory forms  a geometric braid $\beta$, 
a collection of strands in $D^2\times [0,1] $ connecting 
$p_i\times 1$ to $f(p_i) \times 0$ (see Figure~\ref{fig:horseshoe}). 
The isotopy class of~$\beta$ determines the homotopy class of~$f$ 
relative to $\{p_i\} \cup \partial D^2$ and vice versa. 
An $n$-braid refers to the isotopy class of a geometric braid 
with $n$ strands. 
The set of $n$-braids forms the braid group $B_n$.

From now on we consider $f$ as a homeomorphism on a punctured sphere 
$f\colon int(D^2) - \{p_i\} \to int(D^2) - \{p_i\}$. 
In particular by forgetting the boundary $\partial D^2$, 
we ignore Dehn twists along $\partial D^2$ which do not 
affect the dynamics of the braid $\beta$. 
In other words we consider an $n$-braid $\beta$ as 
a mapping class on a $(n+1)$-times punctured sphere 
with the (so called) boundary puncture  fixed. 

\subsection*{Topological entropy}
The topological entropy $h_T(\beta)$ of the braid $\beta$ is 
defined to be the infimum topological entropy of 
the disk homeomorphisms  representing $\beta$. 
$$
h_T(\beta) = \inf_{ g \simeq f} h_T(g) 
$$
The topological entropy of a braid is a conjugacy invariant 
measuring the dynamical complexity of the braid. 
It is equal to the logarithm of the growth rate of the free group 
automorphism induced on $\pi_1(D^2 - \{p_i\})$. 
When $\beta$ is represented by a pseudo-Anosov homeomorphism $f$ 
with dilatation factor $\lambda_f = \lambda(f)$, 
we have $h_T(\beta) = \log \lambda_f$. 
In this case the dilatation $\lambda(\beta)$ of the braid 
is also given by $\lambda(\beta) = \lambda_f$.

If $f$ is homotopic to a periodic homeomorphism, the braid $\beta$ 
is called periodic. 
If there is a collection of disjoint 
sub-surfaces of $int(D^2) - \{p_i\}$ with negative Euler characteristics 
which is homotopically invariant under $f$, 
the braid $\beta$ is called reducible.
As we consider the dynamics of a periodic braid to be trivial, 
studying the dynamics of braids reduces to the maps on 
non-periodic irreducible components.

\subsection*{Pseudo-Anosov homeomorphism}
By the Nielsen-Thurston classification of surface 
homeomorphisms~\cite{MR89k:57023,MR57:16704,MR96d:57014,MR1993745}, 
a non-periodic irreducible braid is represented by a pseudo-Anosov 
homeomorphism. A pseudo-Anosov homeomorphism has several nice 
extremal properties. 
It realizes the minimum topological entropy and the minimum 
quasi-conformality constant in its homotopy class. 
It also has the minimum number of periodic orbits 
for each period~\cite{MR679909}.

A surface homeomorphism $f\colon F\to F$ is called a 
\emph{pseudo-Anosov homeomorphism} 
relative to a puncture set $\{p_i\} \subset F$ when 
the following conditions hold. 
First we need a singular flat metric on $F$ with a finite singularity 
set $\{q_j\}$ such that $\{p_i\} \subset \{q_j\}$. 
Each singularity~$q_j$ has its cone-angle in 
$\{k \pi \mid k\in\mathbb{Z}_{>0}\} $. 
If a singularity has cone-angle $\pi$, it must be one of the 
puncture points~$p_i$. 
Now the homeomorphism $f$ is required to preserve $\{q_j\}$ 
and to be locally  affine (hyperbolic) on $F-\{q_j\}$ 
with a constant dilatation factor $\lambda_f > 1$. 
In particular at a fixed point in $F-\{q_j\}$, 
the map $f$ is locally written  as $\left[ \begin{smallmatrix} 
\lambda_f & 0 \\ 0 & \lambda_f^{-1}
\end{smallmatrix}  \right]$. 

Thus roughly speaking, if a surface homeomorphism $f$ 
represents a non-periodic irreducible mapping class, 
then one can simplify $f$ by pulling tight it everywhere 
until it becomes linear almost everywhere in an appropriate sense.

The horizontal directions to which $f$ expands by the factor $\lambda_f$ 
integrate to form one invariant measured foliation~$\mathcal{F}^s$. 
The vertical directions perpendicular to $\mathcal{F}^s$ 
form the other invariant measured foliation $\mathcal{F}^u$. 
From a singularity $q_j$ with cone-angle $k\pi$, 
$k$-many singular leaves of~$\mathcal{F}^s$ emanate. 
In this case $q_j$ is called a $k$-prong singularity.

Note that in the above definition of pseudo-Anosov homeomorphism 
we can remove or add punctures keeping the same map $f\colon F\to F$. 
When $\{ f^j(x) \} $ is a periodic orbit of non-punctured points,  
\emph{puncturing} at $\{ f^j(x) \}$
refers to adding them to the puncture set~$\{p_i\}$. 
Conversely when $\{ f^j(p_1) \}$ is a periodic orbit of 
$k$-prong punctured singularities  for~$k>1$, 
\emph{capping-off} them refers to removing them from the puncture set. 
For pseudo-Anosov braids, puncturing or capping-off 
corresponds to adding or removing some strands.

Let $\tilde f \colon \tilde F \to \tilde F$ 
be a lift of~$f$ on a finite-fold cover $\tilde F$ of~$F$ 
branched at some finite set of points invariant under~$f$. 
Then by pulling back the flat metric on $F$ to $\tilde F$, 
the lift $\tilde f$ is also a pseudo-Anosov homeomorphism 
with the same dilatation factor $\lambda_{\tilde f} = \lambda_f$.

\subsection*{Train track representative}
\begin{figure}
$$
\epsfig{file=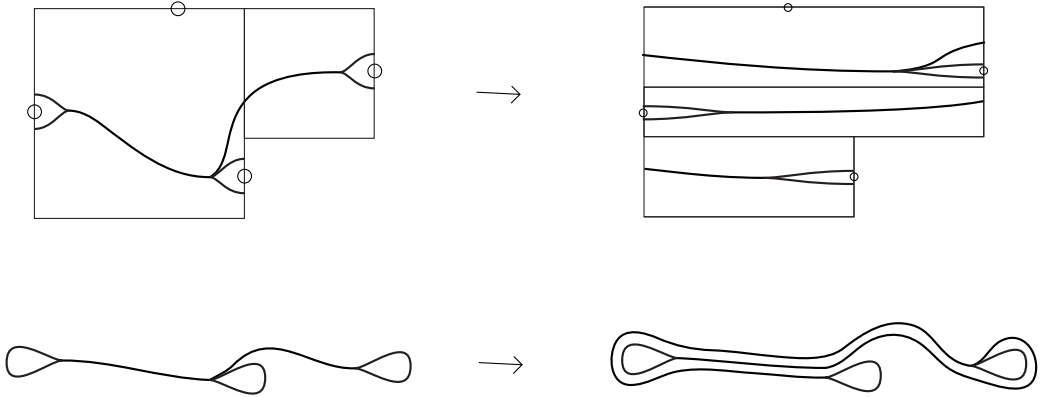,width=0.85\textwidth}
$$
\caption{The boundary of the $L$-shaped region 
reads, counterclockwise from the puncture on the top,  
$abb^{-1}cdd^{-1}eff^{-1}g$ with $ac = g^{-1}e^{-1}$. 
After side-pairing, the $L$-shaped region 
becomes a 4-times punctured sphere with a flat metric.  
The pseudo-Anosov homeomorphism maps 
each rectangle to a longer and thinner horizontal strip 
running over other rectangles.
In the train track representative, 
each edge is assigned tangential and transverse measures 
coming from the width and the height of the corresponding  
rectangle.}\label{fig:markov}

\end{figure}

Using a Markov partition (or its associated train track representative), 
the flat metric and the pseudo-Anosov homeomorphism can be described 
quite concretely (see~\cite[Expos\'e~9]{MR82m:57003} for the definition 
and see Figure~\ref{fig:markov} 
for an example).
Let $\{R_i\}$ be a Markov partition for a pseudo-Anosov homeomorphism $f$. 
The transition matrix $M_f = (m_{ij})$ is defined by 
that $f(R_i)$ crosses over $R_j$ $m_{ij}$-many times. 

The transition matrix $M_f$ is Perron-Frobenius: 
for some $k>0$ each entry of~$M_f^k$ is strictly positive. 
In particular the largest eigenvalue of~$M_f$ is real
and has an eigenvector with strictly positive 
coordinates~\cite[Theorem 1.1 on p.1]{MR0389944}.

The widths $v_i$  and the heights $w_i$ 
of $R_i$ satisfy the following equations.
$$
\lambda_f v_i  = \sum_j m_{ij} v_j  \qquad
w_j  = \sum_i m_{ij} {w_i}/{\lambda_f} 
$$
In particular the dilatation factor $\lambda_f$ appear as 
the eigenvalue of $M_f$ of which  the eigenvector has 
strictly positive coordinates.

We use train track representatives as a notational 
simplification for Markov partitions. 
As in Figure~\ref{fig:markov}, each expanding edge of the invariant 
train track corresponds to a rectangle in the Markov partition. 
Once given the transition matrix of the graph map, 
it is easy to recover the heights and widths of the rectangles.

\subsection*{Main question}
Let us consider the set $\Lambda_{g,n}$ of the dilatation factors 
for pseudo-Anosov homeomorphisms on 
an $n$-times punctured genus-$g$ surface $F_{g,n}$. 
$$
\Lambda_{g,n} = \{\lambda_f \mid f\colon F_{g,n} \to F_{g,n} 
\text{ pseudo-Anosov homeomorphisms}\}
$$
As we can bound the number of rectangles in Markov partitions 
using the Euler characteristic of the punctured surface, 
$\Lambda_{g,n}$ consists of eigenvalues of Perron-Frobenius 
matrices with bounded dimension. 
In particular the set $\Lambda_{g,n}$ is discrete and has 
a minimum. 
Our current work is motivated by the following question. 
\begin{question}\label{quest:main}
What is $\min \Lambda_{g,n}$?
\end{question}
The question asks to find the simplest pseudo-Anosov homeomorphism 
on the surface.

A pseudo-Anosov homeomorphism~$f$ induces an isometry on the 
Teichm\"uller space equipped with the Teichm\"uller metric. 
The pair of invariant measured foliations $(\mathcal{F}^s, \mathcal{F}^u)$ 
determines a geodesic axis in the Teichm\"uller space 
on which $f$ acts as a translation by $\log \lambda_f$.  
The axis projects down to a closed geodesic in the moduli space, 
which is the quotient of the Teichm\"uller space 
by the action of the mapping class group. 
Conversely any closed geodesic in the moduli space 
represents the conjugacy class of some pseudo-Anosov mapping class. 
Therefore Question~\ref{quest:main} can be rephrased 
as asking to find 
the shortest closed geodesic in the moduli space.

The hyperbolic volume of the mapping torus is another 
natural complexity measure for a pseudo-Anosov homeomorphism. 
We notice that a pseudo-Anosov homeomorphism with small dilatation 
tends to have the mapping torus with small hyperbolic volume 
and vice versa.

The question for the minimum volume of the hyperbolic mapping tori 
on a given surface seems to be much more difficult than 
Question~\ref{quest:main}. 
In~\cite{MR1869847} the minimum volume for orientable 
cusped hyperbolic 3-manifolds 
is computed. An extensive use of computer programs is involved 
in its proof. 
In this paper we also use a computer program for the proof 
of the main theorem, but the algorithm and the actual code 
is much simpler than those of~\cite{MR1869847}.

\subsection*{Related results} 
The question asking for the minimum dilatation of pseudo-Anosov 
homeomorphisms on a given surface still remains largely unanswered since after 
the Nielsen-Thurston classification of surface homeomorphisms. 
The existence of Markov partition~\cite[Expos\'e 10]{MR82m:57003} 
for a pseudo-Anosov homeomorphism implies that a dilatation 
should appear as the largest eigenvalue of a Perron-Frobenius 
matrix of bounded dimension, hence in particular should 
be an algebraic integer.  However, it is not clear how the restriction that 
the symbolic dynamical system dictated by a Perron-Frobenius matrix 
is from a homeomorphism on a given surface, actually 
affects the possible values of entropy (the logarithm of the dilatation).

There are several known results relevant to this question of 
the minimum dilatation of pseudo-Anosov homeomorphisms.  
Penner~\cite{MR91m:57010} gives a lower bound  $2^{1/(12 g - 12 + 4n)}$ 
for the dilatations on $F_{g,n}$ a genus-$g$ surface with $n$~punctures. 
In~\cite{MR91m:57010,MR92g:57024,MR2065567} pseudo-Anosov homeomorphisms 
on $F_{g,0}$ with small dilatations are constructed showing 
that the minimum dilatation on $F_{g,0}$ converges to $1$ 
as the genus~$g$ increases. 
Fehrenbach and Los~\cite{MR2001f:37049} compute 
a lower bound $(1+\sqrt2)^{1/n}$ 
for the dilatations of pseudo-Anosov disk homeomorphisms (braids) 
which permute the punctures in one cycle. 
In~\cite{MR2119022} a lower bound $2+\sqrt5$ for the dilatations 
of pseudo-Anosov pure braids is given. 
A pseudo-Anosov disk homeomorphism is represented by a 
transitive Markov tree map preserving the end point set of the tree
with the same topological entropy. 
Baldwin~\cite{MR1838995} gives a lower bound $\log 3$ for the topological 
entropy of transitive Markov tree maps fixing each end point.

The exact values of the minimum dilatations are known only for few 
simple cases. 
Zhirov~\cite{MR1331364} shows that if a pseudo-Anosov homeomorphism 
on $F_{2,0}$ has an orientable invariant foliation, 
then its dilatation is not less than the largest zero $\lambda_5$ of 
$x^4 - x^3 - x^2 - x + 1$, and gives an example of a pseudo-Anosov 
homeomorphism realizing the dilatation~$\lambda_5$.

The pseudo-Anosov 3-braid $\sigma_2\sigma_1^{-1}$ is shown to 
be the minimum in the forcing partial order among 
pseudo-Anosov 3-braids by Matsuoka~\cite{MR860249} 
and Handel~\cite{MR1452182}, 
hence it attains the minimum dilatation. 
The pseudo-Anosov 4-braid $\sigma_3\sigma_2\sigma_1^{-1}$ 
is claimed in~\cite{MR1915500} to have the minimum dilatation, 
but the proof given in~\cite{MR1915500} unfortunately contains an error.

\subsection*{Outline}
In this paper we prove the following theorem 
giving at the same time a corrected proof of the minimality 
of the dilatation of $\sigma_3\sigma_2\sigma_1^{-1} \in B_4$. 

\begin{maintheorem}
The 5-braid 
$\sigma_1\sigma_2\sigma_3\sigma_4\sigma_1\sigma_2$ 
attains the minimum dilatation of pseudo-Anosov 5-braids. 
\end{maintheorem}

The dilatation of a pseudo-Anosov braid is invariant 
under several operations such as conjugation, 
composing with a full twist, taking inverse, and taking reverse. 
It turns out that for braid indices 3 to 5, 
the pseudo-Anosov braids realizing the minimum dilatations 
are essentially unique, modulo the aforementioned operations. 
This could be just a coincidence. 
It would be a good surprise if some uniqueness property 
can be proved for the minimum-dilatation pseudo-Anosov braids.

The two main ingredients of the proof of Theorem~\ref{thm:main} 
are the construction of folding automata for generating 
candidate pseudo-Anosov braids for the minimum dilatation, 
and the following lemma for bounding the word lengths of the 
candidate braids.

\begin{mainlemma}
If $M$ is a Perron-Frobenius matrix of dimension $n$ 
with  $\lambda>1$ its largest eigenvalue, 
then 
$$
\lambda^n \ge |M| - n + 1
$$
where $|M| $ denotes the sum of entries of $M$. 
\end{mainlemma}
This lemma improves on \cite[Theorem~6]{MR1038734} 
and replaces erroneous Lemma~3,4 of~\cite{MR1915500}.

Given a pseudo-Anosov homeomorphism $f\colon (F,\{p_i\}) \to (F,\{p_i\})$ 
on a surface $F$ with punctures $p_i$ with 
negative Euler characteristic $\chi(F-\{ p_i \}) < 0$, 
there exists a  train track representative of $f$. 
There exists an invariant train track $\tau \subset F - \{ p_i \}$ 
which carries $f(\tau)$.  
In particular there is a splitting sequence 
$$
\tau  = \tau_0 \succ \tau_1 \succ \cdots \succ \tau_k  = f(\tau)
$$
from $\tau$ to $f(\tau)$, 
where $\tau_j \succ \tau_{j+1}$ is an elementary splitting move. 

By observing that there are only finitely many diffeomorphism types 
of the pair $( F - \{p_i\}, \tau_j)$, 
one can effectively construct a \emph{splitting automaton}, 
which is a finite graph with train tracks as its vertices 
and with splitting moves as its arrows.

The existence of the train track representative, in particular 
of the splitting sequence, implies that every pseudo-Anosov homeomorphism 
appears, up to conjugacy, as a closed path in some splitting 
automata (see~\cite{MR914107}).  
Papadopoulos and Penner~\cite[Theorem~6]{MR1038734} also gave a lower bound 
for the dilatation in terms of word length in automata. 

In this paper we use folding automata as in~\cite{MR1915500}, 
which are finite graphs with embedded train tracks as vertices 
and with elementary folding maps as arrows. 
An elementary folding map is an inverse of a splitting move.

If  we are given an upper  bound for the word length
in terms of the dilatation, then 
on a fixed folding automaton, the search for the minimum dilatation 
in the automaton reduces to checking for finitely many 
closed paths.

Lemma~\ref{lemma:bound}, which is an improvement of 
\cite[Theorem~6]{MR1038734}, not only gives an 
upper bound of the word lengths of mapping classes with dilatation 
bounded by a fixed number, 
but also trims out many branches which appear in the course of search 
in a big tree, namely the set of paths with bounded length.  
In fact Lemma~\ref{lemma:bound} implies  that it suffices to 
consider only such paths whose any subpath has a transition matrix 
with bounded norm.  

For the minimum dilatation of 5-braids, the previously mentioned 
restriction on paths by transition matrix norm  and 
another restriction by Lemma~\ref{lem:avoidword} 
significantly reduce the number of candidate 
braids  making the computation feasible. 

We think that the same method for computing the minimum dilatation 
would still work for a few more simple cases like  on a genus-2 closed 
surface, 
although it would involve more complicated computer aided search.

\section{Folding automata}
Given a pseudo-Anosov homeomorphism $f\colon (F, \{p_i\}) \to (F, \{p_i\})$
on a closed surface $F$ with punctures $\{p_i\}$, 
there exists an invariant train track $\tau \subset F - \{p_i\}$ 
and $f$ is represented by a train track 
map $f_\tau \colon \tau \to \tau$~\cite{MR914107}. 

\begin{figure}
$$
\epsfig{file=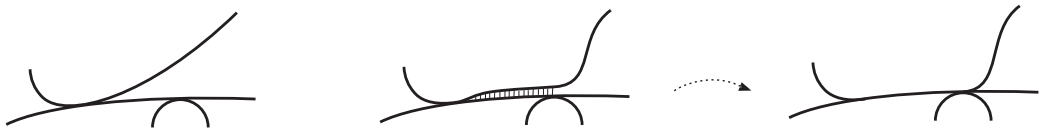,width=0.85\textwidth}
$$
\caption{An elementary folding map}\label{fig:folding}
\end{figure}

A train track $\tau$ is a smooth branched 1-manifold embedded in 
the surface $F-\{p_i\}$ such that each component of the complement 
$F - \{p_i\} - \tau$ is 
either a once punctured $k$-gon for~$k\ge 1$ or a 
non-punctured $k$-gon for~$k\ge 3$. 
The train track $\tau$ is called \emph{invariant} 
under~$f$ if~$f(\tau)$ smoothly collapses onto $\tau$ in $F-\{p_i\}$ 
inducing a smooth map $f_\tau \colon \tau \to \tau$, 
which maps branch points to branch points. 
In this case 
one may repeatedly fold (or zip) $f(\tau)$ nearby cusps to obtain 
a train track isotopic to $\tau$ in $F - \{p_i\}$ 
(see Figure~\ref{fig:folding} and \cite[Fig.~4,~5]{MR1628757}). 

Let $f_\tau \colon \tau \to \tau$ be a train track representative 
of a pseudo-Anosov homeomorphism $f$. 
An edge $e$ of $\tau$ is called \emph{infinitesimal} 
if it is eventually periodic under $f_\tau$, that is,  
$f_\tau^{N+k} (e) = f_\tau^N (e)$ for some $N,k>0$. 
An edge  of $\tau$ is called \emph{expanding} if it is not 
infinitesimal. 
An expanding edge $e$ actually has a positive length in the sense 
that $\lim_{N\to\infty} |f_\tau^N(e) | / \lambda_f^N$ is positive 
where $|\cdot|$ denotes the word length of a path 
and $\lambda_f = \lambda(f)$ denotes the dilatation factor for $f$.

A graph map is called \emph{Markov} if it maps vertices to 
vertices, and is locally injective at points that do not map into 
vertices. 
Given a Markov map $g\colon \tau \to \tau'$, 
the transition matrix $M_g = (m_{ij})$ is 
defined by that the $j$-th edge $(e'_j)^{\pm1}$ of $\tau'$ 
occurs $m_{ij}$-many times in the path $g(e_i)$, 
the image of the $i$-th edge of~$\tau$.
When $\tau' = \tau$, the transition matrix is square and 
considering its spectral radius makes sense. 

The spectral radius of $M_{f_\tau}$ equals the dilatations 
factor $\lambda(f)$ for the pseudo-Anosov homeomorphism $f$. 
Coordinates of the corresponding eigenvectors of 
$M_{f_\tau}$ and its transpose $M_{f_\tau}^T$ are 
tangential and transverse measures of edges of $\tau$, 
which are projectively invariant under~$f$.

An elementary folding map $\pi\colon \tau \to \tau'$ 
is a smooth Markov map between two train tracks $\tau$ and $\tau'$
such that for only one edge $e$ of $\tau$, the image $\pi(e)$ 
has word length $2$, and the other edges map to paths of length $1$. 
In other words the transition matrix $M_\pi$ is of the form 
$P+B$ for some permutation matrix $P$ and for 
some elementary matrix $B$.

When the train tracks are embedded in a surface as in our case of concern, 
the pairs of edges which are folded should be adjacent 
in the surface: the two segments of $\tau$ which are identified 
by the elementary folding map are two sides of an open triangle 
in $F - \{p_i\} - \tau$ (see Figure~\ref{fig:folding}).

\begin{prop}
A train track representative $f_\tau \colon \tau \to \tau$ of 
a surface homeomorphism $f\colon (F , \{p_i\}) \to (F , \{p_i\})$ 
admits a folding decomposition as follows: 
$$
f_\tau = \rho \circ \pi_k \circ \cdots \circ \pi_1
$$
where $\pi_j \colon \tau_j \to \tau_{j+1}$ are 
elementary folding maps, 
$\tau_1  = \tau_{k+1} = \tau$, 
and $\rho\colon \tau \to \tau$ is an isomorphism 
induced by a periodic surface homeomorphism $(F - \{p_i\}, \tau) \to 
(F - \{p_i\}, \tau)$ preserving $\tau$. 
\end{prop}
\begin{proof}
It follows from~\cite{MR85m:05037a}. 
See~\cite{MR914107,MR1915500} for more details. 
\end{proof}

By observing that there are only finitely many 
possible diffeomorphism types for the pairs $(F-\{p_i\},  \tau_j)$ 
appearing in the folding decomposition, 
we can construct folding automata. 
A \emph{folding automaton} is a connected directed graph 
with diffeomorphism types of train tracks as vertices, 
with elementary folding maps and isomorphisms as arrows. 
See Figure~\ref{fig:twotrigon} for a simplified version 
of a folding automaton. 
The train tracks in Figure~\ref{fig:twotrigon} admit no 
non-trivial isomorphisms, that is, 
if $h\colon (D_5,\tau) \to (D_5,\tau)$ is an 
orientation preserving diffeomorphism fixing~$\tau$ in the automaton, 
then $h$ is isotopic to the identity map.  
So in Figure~\ref{fig:twotrigon} there are no arrows 
corresponding to isomorphisms.

\begin{coro}
All train track representatives of  pseudo-Anosov homeomorphisms 
are represented by  closed oriented paths in folding automata. 
\end{coro}

To each closed path based at a train track $\tau$ in  a folding 
automaton, 
associated is a train track representative $f_\tau \colon \tau \to \tau$ 
of some homeomorphism $f\colon (F , \{p_i\}) \to (F , \{p_i\})$. 
The disk homeomorphism $f$ is pseudo-Anosov if and only if 
the transition matrix $M_{f_\tau}$ is Perron-Frobenius 
(also called primitive) modulo infinitesimal edges: for some $N>0$, 
the power $M_{f_\tau}^N$  is strictly positive in the block of 
expanding edges. 
To find out whether $M$ is Perron-Frobenius, 
it suffices by~\cite{MR0097416,MR0035265}, 
\cite[Theorem 2.8 on p.52]{MR0389944}  to check 
if $M^{n^2 - 2n + 2}$ has all non-zero entries where $n$ is the dimension 
of the matrix $M$.

Now we discuss simplifying the train track maps so that 
we can restrict to simplified folding automata. 
If the pseudo-Anosov homeomorphism $f$ fixes a distinguished 
puncture~$p_0$: $f(p_0) = p_0$ (for instance when $f$ is from a disk 
homeomorphism and $p_0$ is the boundary puncture), 
we can give a restriction to the train track map 
$f_\tau \colon \tau \to \tau$, thereby reducing 
the size of the folding automata needed in our computation.

We first assume that only the component of $F-\tau$ containing $p_0$ 
has expanding edges on its sides: 
the other components of $F-\tau$ not containing $p_0$ 
are bounded only by infinitesimal edges. 
If one is given a train track representative $f_\tau \colon \tau \to \tau$ 
not satisfying this assumption, 
he can apply a splitting operation~\cite[Section~5]{MR96d:57014} nearby 
$p_0$ (when $p_0$ is enclosed only by infinitesimal edges) 
then apply a sequence of folding 
operations~\cite[p.15]{MR92m:20017} 
\cite[Section~2.2]{MR97h:20031} nearby other punctures $p_i$, $i\neq 0$, 
until all the components of the train track complement  not 
containing $p_0$ shrink to be infinitesimal, 
to obtain a new train track representative satisfying 
the assumption~\cite[Proposition~3.3]{MR92m:20017}.

\begin{figure}
$$
\epsfig{file=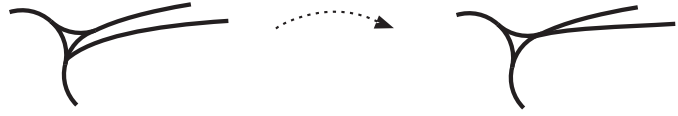,width=0.55\textwidth}
$$
\caption{A folding operation which absorbs 
an infinitesimal edge}\label{fig:foldinfini}
\end{figure}

Applying some more folding operations~(see Figure~\ref{fig:foldinfini}), 
we can also remove 
any cusp occurring between an expanding edge and an infinitesimal edge.

\begin{figure}
$$
\epsfig{file=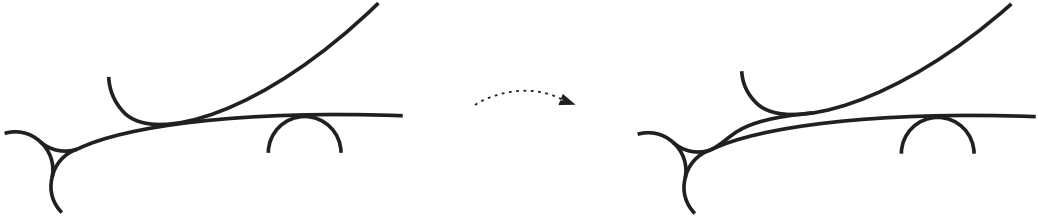,width=0.85\textwidth}
$$
\caption{A splitting operation which shifts a cusp 
incident only to expanding edges}\label{fig:splitting}
\end{figure}

We assume that  cusps  occur only at corners of infinitesimal 
multigons. 
If one is given a train track representative 
with a cusp  incident only to expanding edges,
not satisfying this assumption, 
he can apply a splitting operation at the cusp  until 
the cusp hits an infinitesimal multigon (see Figure~\ref{fig:splitting}).

\begin{figure}
$$
\epsfig{file=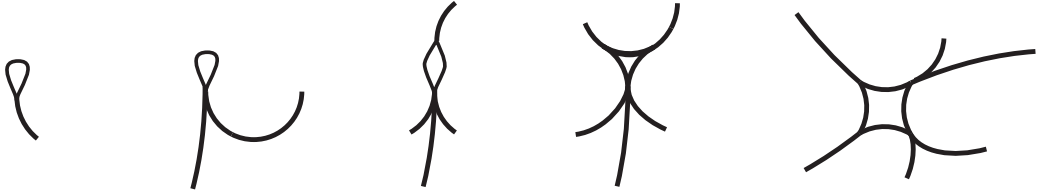,width=0.85\textwidth}
$$
\caption{Local models for train tracks in simplified folding 
automata}\label{fig:legmultigon}
\end{figure}

\begin{figure}
$$
\epsfig{file=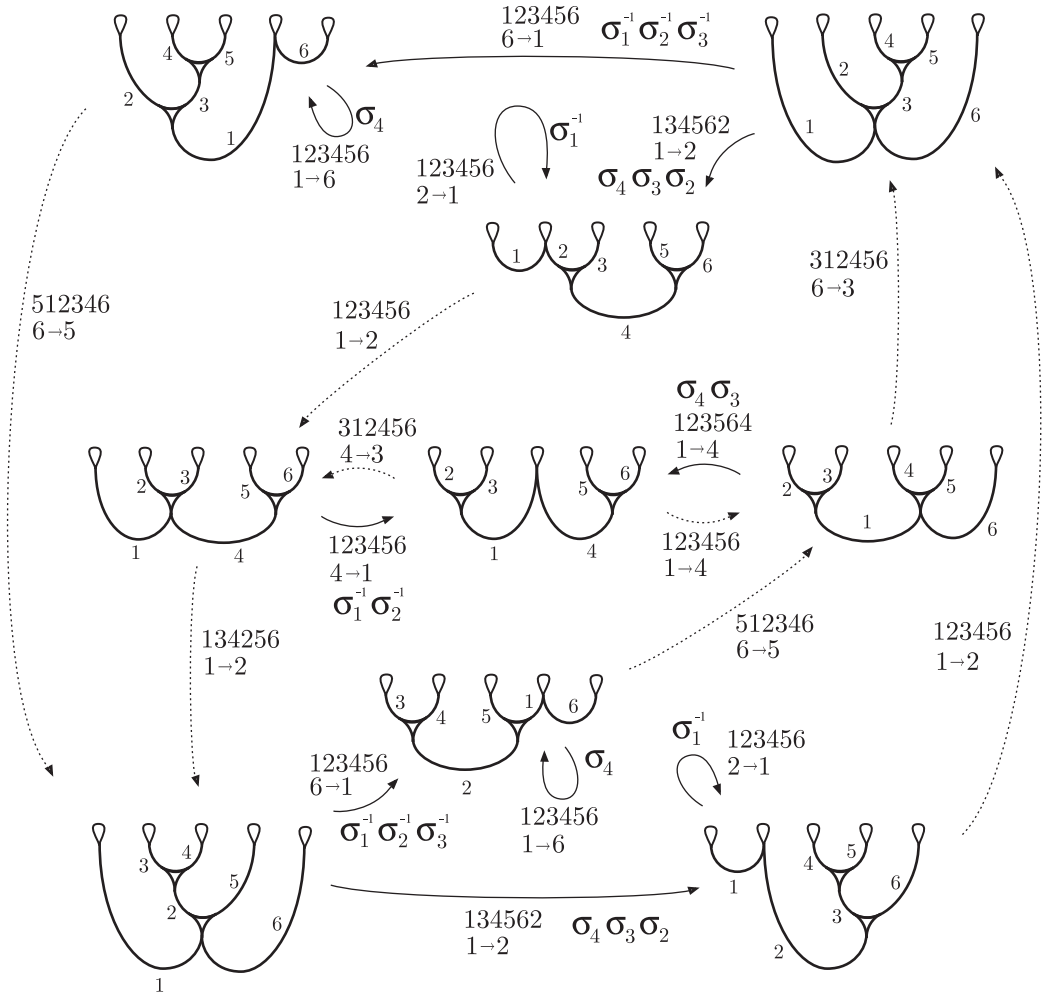,width=0.85\textwidth}
$$
\caption{A folding automaton for pseudo-Anosov 5-braids 
with two non-punctured 3-prong singularities}\label{fig:twotrigon}
\end{figure}

Therefore a pseudo-Anosov braid has an invariant train track 
which is locally modeled by 
infinitesimal $k$-gons to which expanding edges are joined 
(possibly) forming cusps only between expanding edges 
(see Figure~\ref{fig:legmultigon}).

In this paper we use simplified versions of folding automata, 
of which train tracks satisfy the previously given conditions,  
and each arrow is either an isomorphism or 
a composite of two elementary folding maps 
by which one expanding edge and one infinitesimal edge is absorbed 
into another expanding edge. 
It is not hard to see that simplified folding automata also 
generate all the conjugacy classes of pseudo-Anosov homeomorphisms.

In this paper our subject of interest is pseudo-Anosov homeomorphisms 
on a 5-times punctured disk $D_5$, 
or equivalently on a 6-times punctured sphere $F_{0,6}$ with a 
distinguished boundary puncture.

We explain how to  read Figure~\ref{fig:twotrigon}, which 
depicts a simplified version of a folding automaton. 
Each train track is embedded in a 5-times punctured disk, 
with each puncture enclosed by an infinitesimal monogon. 
Each embedding is chosen arbitrarily, and only  
the orientation preserving diffeomorphism types of embedded 
train tracks count. 

An arrow is a composite of two elementary folding maps, 
one involving an infinitesimal edge and another involving  
only expanding edges. 
We ignore the infinitesimal edges in computing the transition matrix 
because the occurrences of infinitesimal edges do not affect 
the resulting dilatation factor. 

An arrow is drawn dashed if it induces a homeomorphism 
isotopic to identity, and it is drawn by a solid line otherwise. 
Note that a folding map $\pi\colon \tau \to \tau'$ determines 
a disk homeomorphism $f\colon D_5 \to D_5$ up to isotopy 
when the embeddings $\tau \to D_5$ and $\tau' \to D_5$ of 
the two train tracks  are fixed.  
In particular $\tau' \succ f(\tau)$, that is, 
$f(\tau)$ folds to be $\tau'$ inducing the folding map $\pi$.
To each solid arrow, a braid word is assigned representing 
the associated disk homeomorphism. 

Edges of a train track are numbered by $\{1,2,\cdots,6\}$ 
in such a way that  
in the peripheral word running clockwise from a cusp, 
new edges appear in an increasing order.  
This naming of edges amounts to fixing a groupoid homomorphism 
from paths in the automaton to transition matrices, 
that is, for two paths $\gamma$ and $\delta$, 
$M(\gamma\cdot\delta) = M(\gamma)M(\delta)$ if $\gamma$ ends 
at the starting vertex of $\delta$, where $M(\gamma)$ denotes 
the transition matrix for $\gamma$.

Each arrow is associated with a permutation 
$i_1 i_2 i_3 i_4 i_5 i_6 $ and a rule $m\to n$, meaning  that
under the elementary folding map, the edge $j$ maps to~$i_j$ for 
$j\neq m$, and $m$ maps to~$i_m\cdot n$. 
(Here we concern only the transition matrix so that 
the direction of edges and the order of concatenation 
are irrelevant.)

Given an adjacent pair $(e_1, e_2)$ of edges with a cusp 
between them, there are two possible folding maps: 
one under which the image of $e_1$ overpasses that of $e_2$ 
and another vice versa.  
Therefore from each train track in Figure~\ref{fig:twotrigon}, 
two arrows of elementary folding maps emanate. 
Likewise two arrows are headed for each train track, 
because at each cusp there are two different elementary splittings 
possible.

\section{Search for the minimum dilatation} 
In this section we prove that the largest zero $\lambda_5$ 
of $x^4 - x^3 - x^2 - x + 1$ is indeed the minimum dilatation 
for pseudo-Anosov 5-braids. 

The problem for the minimum dilatation reduces to a search 
in a finite set of closed paths in folding automata 
because by~\cite[Theorem~6]{MR1038734} or by Lemma~\ref{lemma:bound}  
the dilatation grows as the norm of the transition matrix 
grows, 
and there are only finitely many closed paths whose transition matrices
have norm bounded by a given number. 
For instance if a closed path in folding automata has length $N$, 
then its associated transition matrix has norm at least $N$.

We first restate and prove the lemma in Section~1. 
\begin{lemma}\label{lemma:bound}
If $M$ is a Perron-Frobenius matrix of dimension $n$ 
with  $\lambda>1$ its largest eigenvalue, 
then 
$$
\lambda^n \ge |M| - n + 1
$$
where $|M| $ denotes the sum of entries of $M$. 
\end{lemma}
\begin{proof}
Let $M = (m_{ij})$ and  let
$(v_i)$ be the eigenvector given by the equation 
$$\lambda v_i = \sum_{j=1}^n m_{ij} v_j$$ 
for $v_i>0$, $1\le i \le n$.

The matrix $M$ is the transition matrix of a graph $G$ 
with vertex set $V(G) = \{1,2,\ldots, n\}$  such that 
the number of oriented edges from $i$ to $j$ equals $m_{ij}$. 
Let $M^n = (k_{ij})$.  The number of paths with length $n$ 
from $i$ to $j$, equals $k_{ij}$. 
For each pair $(i,j)$ of vertices there exists an oriented path 
from $i$ to $j$ since $M$ is Perron-Frobenius.

Note that $(v_i)$ is also the eigenvector of $M^n$ 
with eigenvalue $\lambda^n$.  
Choose  $v_p = \min_i v_i$ the smallest coordinate of $(v_i)$. 
\begin{align*}
\lambda^n v_p &= \sum_j k_{pj} v_j   \\
&\ge (\sum_j k_{pj}) v_p \qquad \text{since $v_j \ge v_p$}  
\end{align*}
The inequality $\lambda^n \ge \sum_j k_{pj}$ reads 
that $\lambda^n$ is not less than the number of length-$n$ paths from 
the vertex $p$ of $G$. 

Take a maximal positive tree $T\subset G$ rooted at $p${\;}: 
 each vertex of $G$ is connected  to $p$  by a 
unique oriented path in~$T$. 
As $T$ is maximal, $|V(T)| =  |V(G)| = n $ 
so that the number of edges of $T$ is $|E(T)| = n - 1$. 

If an oriented path $\gamma$ from $p$ has length $n$, 
it must digress from $T$ at some point. 
Define $a(\gamma) = e \in E(G) - E(T)$ to be 
the first edge \emph{not} in $E(T)$ of $\gamma$.  
This defines a function 
$$
a \colon \{ \text{length-$n$ paths from $p$ in $G$} \}  \to 
E(G) - E(T). 
$$
Clearly $a$ is a surjection since the tail of each edge in $E(G) - E(T)$ 
is connected to $p$  by a path in $T$ with length at most $n-1$. 
Therefore 
$$
\lambda^n \ge \sum_j k_{pj} \ge | E(G) - E(T) | =  |M| - (n-1). 
$$
\end{proof}

\begin{thm}[\cite{MR1915500}]\label{thm:fourbraid}
If a pseudo-Anosov 4-braid has an invariant foliation with 
one non-punctured 3-prong singularity, then 
its dilatation is not less than $\lambda_4 \approx 2.29663$, 
the largest zero of $x^4 - 2 x^3 - 2x + 1$. 
\end{thm}
\begin{proof}
Let $\beta \in B_4$ be a pseudo-Anosov 4-braid 
in the theorem.  
Up to conjugacy and multiplication by central elements, 
$\beta$ appears as a closed path $\gamma$ in the folding 
automaton in Figure~\ref{fig:onetrigon}.

\begin{figure}
$$
\epsfig{file=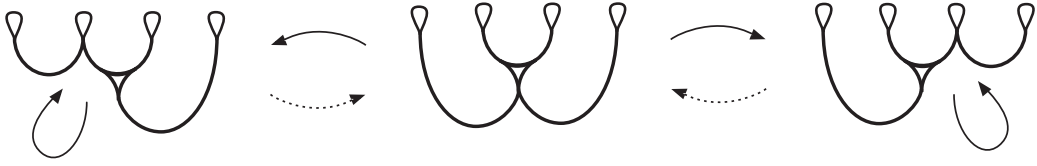,width=0.85\textwidth}
$$
\caption{A simplified folding automaton for pseudo-Anosov 4-braids with 
a 3-prong singularity}\label{fig:onetrigon}
\end{figure}

If $\lambda(\beta) < 2.3$, 
then by Lemma~\ref{lemma:bound} 
we have a bound $31 > 2.3^4 + 4 - 1 \ge |M_\gamma|$ 
for the norm of its transition matrix $M_\gamma$. 

By a computer aided search~\cite{fbrmin} 
in the finite set of paths $\gamma$ with $|M_\gamma| < 31$, 
we conclude that up to conjugacy, multiplication by 
central elements, and taking inverse, 
the braid 
$\sigma_3\sigma_2 \sigma_1^{-1}$ is the only pseudo-Anosov  
 4-braid  with dilatation less than~$2.3$.
It can be easily checked that $\lambda(\sigma_3\sigma_2 \sigma_1^{-1}) = 
\lambda_4$. 
\end{proof}

We say two matrices have the same \emph{pattern} 
if they have zero entries and positive entries in 
the same positions.  
We write $M \ge M'$ for $M = (m_{ij})$ and $M' = (m'_{ij})$ 
if $m_{ij} \ge m'_{ij}$ for all $i,j$. 

In the following lemma we ignore the parts of transition 
matrices arising from infinitesimal edges, 
so that for a closed path to 
represent a pseudo-Anosov  homeomorphism implies   
for its transition matrix to be Perron-Frobenius. 
\begin{lemma}\label{lem:avoidword}
Let $\gamma$ be a closed path  in a folding automaton. 
Let $N, k>0$ be numbers such that 
the transition matrices 
$M(\gamma^{N + i + k})$ and $M(\gamma^{N + i})$  have 
the same pattern  and $M(\gamma^{N + i + k}) \ge M(\gamma^{N + i})$  
for any $i\ge 0$. 
Then a closed path of the form $\alpha \cdot \gamma^{N+i + k} \cdot 
\delta$ in the folding automaton 
represents a pseudo-Anosov homeomorphism  
if and only if $\alpha \cdot \gamma^{N+i } \cdot \delta$  does. 
Furthermore  in this case 
we have an inequality 
$$ \lambda(\alpha \cdot \gamma^{N+i } \cdot \delta) \le 
 \lambda(\alpha \cdot \gamma^{N+i + k } \cdot \delta)  $$
between their dilatation factors.
\end{lemma}
\begin{proof}
It suffices to prove the lemma  for 
$  \delta \cdot\alpha \cdot \gamma^{N+i + k } $ and 
$  \delta \cdot\alpha \cdot \gamma^{N+i } $
since conjugation does not affect dilatation factor or being pseudo-Anosov.

Since $M(\gamma^{N+i + k })$ and $M(\gamma^{N+i })$ 
have the same pattern, 
$M(\delta\cdot\alpha)  M(\gamma^{N+i + k })$ and 
$M(\delta\cdot\alpha)  M(\gamma^{N+i })$  also
have the same pattern. 
In particular one is Perron-Frobenius if and only if so is the other, 
which proves the first claim of the lemma.

From $  M(\gamma^{N+i  }) \le   M(\gamma^{N+i+k}) $, 
we have 
 $$  M(\delta \cdot\alpha \cdot \gamma^{N+i  }) \le  
  M(\delta \cdot\alpha \cdot \gamma^{N+i + k }), $$ 
which by~\cite[Theorem 1.1 (e)]{MR0389944} 
implies the inequality of the lemma.
\end{proof}

\begin{remark}\label{rem:trim}
Let $\gamma, N, k$ be given as in Lemma~\ref{lem:avoidword}. 
Then the lemma implies that 
when we search just for  the minimum dilatation factor for 
pseudo-Anosov homeomorphisms, 
it suffices to search in the set of paths 
that do not contain $\gamma^{N+k}$ as a subpath. 

In the search in the automaton in Figure~\ref{fig:twotrigon}, 
we exclude paths containing several closed paths, 
for example $ \left({{123564}\atop{1\to4~~~~}} \cdot  
{{123456}\atop{1\to4~~~~}}\right)^6 $, 
$ \left({{123456}\atop{4\to1~~~~}} \cdot  
{{312456}\atop{4\to3~~~~}}\right)^6 $, 
and second iterates of  length 1 loops. 
This reduces the size of the set of candidate braids for 
minimum dilatation to the extent that 
the computation in the proof of Theorem~\ref{thm:twotrigon} 
becomes possible on a personal computer. 
\end{remark}

\begin{thm}\label{thm:twotrigon}
If a pseudo-Anosov 5-braid has an invariant foliation with 
two non-punctured 3-prong singularities, then 
its dilatation is not less than 
the largest zero 
$\approx 2.01536$ of $x^6 -  x^5 -  4 x^3  - x + 1$. 
\end{thm}
\begin{proof}
It is easy to check that there are only nine 
different diffeomorphism types of train tracks in $D_5$, 
locally modeled by infinitesimal multigons 
with outgoing expanding-edge legs as in Figure~\ref{fig:legmultigon}. 
By computing the elementary folding maps among them 
(more precisely composites of two elementary folding maps, 
one of them involving an infinitesimal edge), 
we have a folding automaton depicted in Figure~\ref{fig:twotrigon}. 

By a computer aided search~\cite{fbrmin} in the set of 
paths $\gamma$ with $|M_\gamma| \le 2.02^6 + 5 < 73$, 
we conclude that up to 
conjugacy and multiplication by central elements, 
$\sigma_4\sigma_3\sigma_1^{-1}\sigma_2^{-1}$ with 
dilatation $\approx 2.01536$ is the only 
such pseudo-Anosov 5-braid with dilatation less than~$2.02$. 
\end{proof}

\begin{lemma}\label{lem:fourgon}
If a pseudo-Anosov 5-braid has an invariant foliation with 
a non-punctured 4-prong singularity, then 
its dilatation is not less than  
the largest zero $\approx 2.15372$  
of $x^4 -  3 x^3  + 3 x^2   -  3 x + 1$. 
\end{lemma}
\begin{proof}
The folding automaton  for this case is similar to 
one in Figure~\ref{fig:onetrigon}. 
As in the proof of Theorem~\ref{thm:fourbraid}, 
a computer aided search~\cite{fbrmin} in the set of 
closed paths up to length $56 > 2.2^5 + (5-1)$ 
shows that the largest zero of $x^4 -  3 x^3  + 3 x^2   -  3 x + 1$
is the minimum dilatation factor in the automaton, 
and it is achieved by $\sigma_4\sigma_3\sigma_2\sigma_1^{-1}$. 
\end{proof}

\begin{lemma}\label{lem:punctrigon}
If a pseudo-Anosov 5-braid has an invariant foliation with 
a punctured 3-prong singularity, then 
its dilatation is not less than  $\lambda_5$, 
the largest zero  of $x^4 -  x^3 -  x^2   - x + 1$. 
\end{lemma}
\begin{proof}

Let $f\colon F_{0,6} \to F_{0,6}$ be a pseudo-Anosov homeomorphism 
with an invariant foliation~$\mathcal{F}$. 
If the invariant measured foliation~$\mathcal{F}$ 
on a 6-times punctured sphere
has a punctured 3-prong singularity, 
then it has five other punctured 1-prong singularities and no more.

We can assume that the punctured 3-prong singularity is  
the boundary puncture since it should be fixed by the homeomorphism $f$. 
Now we use the folding automaton  that generates such pseudo-Anosov 
braids with  three prongs at the boundary puncture.

There are eleven diffeomorphism types of train tracks to consider 
for this case (see Figure~\ref{fig:threecusp}). 
There are fifty arrows in the automaton, 
which are too many to be drawn in a figure in this paper. 
See~\cite{fbrmin} for details.

\begin{figure}
$$
\epsfig{file=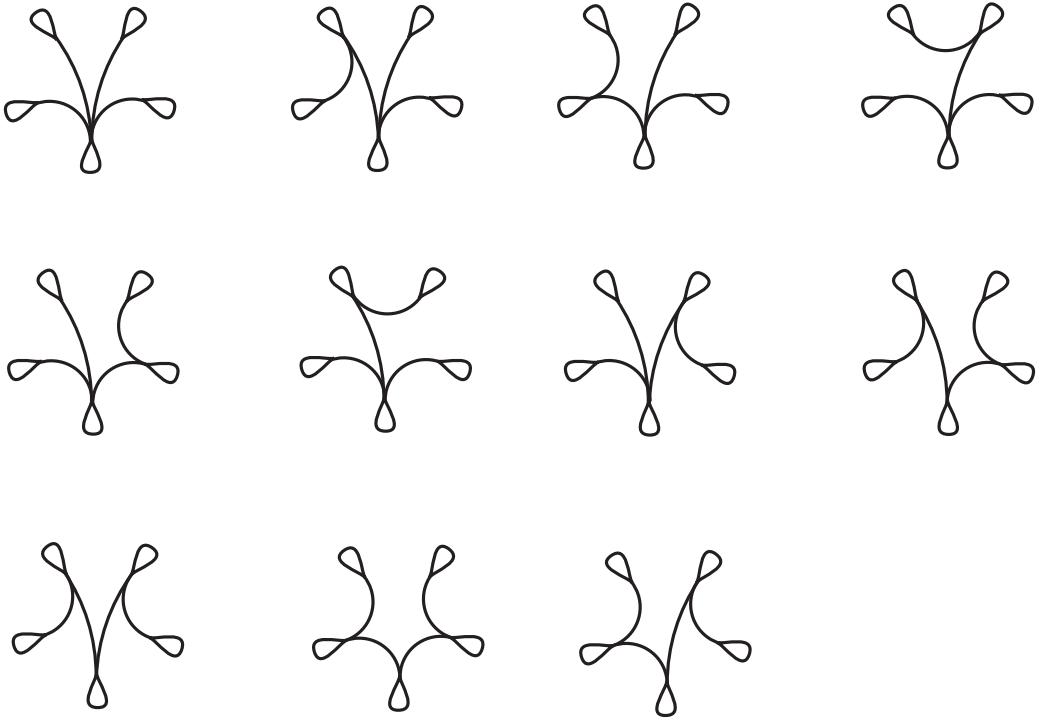,width=0.85\textwidth}
$$
\caption{Train tracks for pseudo-Anosov 5-braids with a 3-pronged 
boundary puncture}\label{fig:threecusp}
\end{figure}

By the same kind of computer aided search as before, 
in the set of closed paths in the folding automaton 
up to length  $12 \ge (\lambda_5)^4 + (4-1)$, 
we conclude that $\lambda_5$ is the minimum dilatation 
factor for pseudo-Anosov braids in this automaton. 

The dilatation is achieved by 
$\sigma_1\sigma_2\sigma_3\sigma_4\sigma_1\sigma_2$.
\end{proof}

\begin{remark}\label{rem:zhirov}
In~\cite{MR1331364}, $\lambda_5$ is proved to be 
the minimum dilatation factor for a pseudo-Anosov homeomorphism 
with an orientable invariant foliation on a closed 
genus-2 surface. 
The proof in~\cite{MR1331364} seems to have a gap. 
To complete the proof one needs to show that 
the golden ratio $(1+\sqrt5)/2 \approx 1.61803$  
the largest zero of $x^4 - 3x^2 + 1$, cannot be a dilatation 
factor for such a pseudo-Anosov homeomorphism. 
In~\cite{MR1485474}, it is proved that 
such a pseudo-Anosov homeomorphism with quadratic dilatation factor 
is a lift of an Anosov homeomorphism via a branched covering. 
Lemma~\ref{lem:fourgon} and Lemma~\ref{lem:punctrigon} 
follow from~\cite{MR1331364} 
by taking double covers branched at odd-prong singularities. 
\end{remark}

By collecting all the results, 
we conclude this paper with a proof of the main theorem.

\begin{thm}\label{thm:main}
The 5-braid 
$\sigma_1\sigma_2\sigma_3\sigma_4\sigma_1\sigma_2$ 
attains the minimum dilatation of pseudo-Anosov 5-braids. 
\end{thm} 
\begin{proof}
Let $f\colon F_{0,6} \to F_{0,6}$ be a pseudo-Anosov 
homeomorphism on a 6-times punctured sphere 
with a punctured point fixed by $f$. 
Let~$\mathcal{F}$ be its invariant measured foliation. 
Since~$\mathcal{F}$ has exactly six punctures, 
the formula $2 = \chi(F_{0,0}) = \sum_k (1 - k/2) n_k$ 
where $n_k$ denotes the number of $k$-prong singularities, 
says no singularity of~$\mathcal{F}$ can have more than four prongs. 

The list of possible types of~$\mathcal{F}$ according to its 
singularity type is: 
\begin{enumerate}
\item six punctured 1-prong singularities and one non-punctured 
4-prong singularity: $n_1 = 6,\ n_4=1$ \label{item1}
\item six punctured 1-prong singularities and two non-punctured 
3-prong singularities: $n_1 = 6,\ n_3=2$ \label{item2}
\item five punctured 1-prong singularities and one punctured 
3-prong singularity: $n_1 = 5,\ n_3=1$ \label{item3}
\item five punctured 1-prong singularities,  one punctured 2-prong 
singularity, and one non-punctured 
3-prong singularities: $n_1 = 5,\ n_2=1,\ n_3=1$ \label{item5}
\item four punctured 1-prong singularities and two punctured 
2-prong singularities: $n_1 = 4,\ n_2=2$ \label{item4}
\end{enumerate}

The case~(\ref{item4})  is, by capping-off 2-prong singularity punctures, 
of a pseudo-Anosov homeomorphism on a four-times punctured sphere, 
which lifts to an Anosov homeomorphism on a torus via 
branched double covering. 
Therefore in this case $\lambda(f) \ge (3 + \sqrt5)/2 > \lambda_5$.

The case~(\ref{item5}) reduces to the case~(\ref{item3}) 
by capping-off the punctured 2-prong singularity and 
puncturing at the 3-prong singularity.

For the cases (\ref{item1}) and (\ref{item3}) 
we have $\lambda(f) \ge \lambda_5$ by 
Lemma~\ref{lem:fourgon} and Lemma~\ref{lem:punctrigon} 
or by Remark~\ref{rem:zhirov}.

Finally the case~(\ref{item2}) is covered by Theorem~\ref{thm:twotrigon}, 
so that we have $\lambda(f) > 2.01 > \lambda_5$. 
In fact this is the only part of the proof which actually 
requires a computer aided search 
if one uses Zhirov's result~\cite{MR1331364}.

Collecting all of these we conclude that 
$\lambda(f) \ge \lambda_5 \approx 1.72208$. 
It is easily checked that $\beta =  
\sigma_1\sigma_2\sigma_3\sigma_4\sigma_1\sigma_2$ 
realizes this dilatation $\lambda(\beta) = \lambda_5$. 
\end{proof}

\section{Implementation} 
To search for the minimum-dilatation pseudo-Anosov homeomorphism 
on a given surface, we first need to generate a collection 
of folding automata. For 5-braids it is possible to build the necessary 
folding automata manually. On surfaces with more punctures and greater 
genus, we also need a computer program \texttt{genauto} to generate 
the folding automata. This paper will not cover the details 
of its implementation. The following is a pseudo-code for \texttt{genauto}.

\begin{center}
\begin{boxedminipage}{0.8\textwidth}
\vspace*{0.5\topskip}
\begin{center}\texttt{genauto}\end{center}
\hskip1.5em
\begin{minipage}{0.85\textwidth}
\begin{itemize}
\item[\textbf{input}:] genus $g$ and the number of punctures $n$
\item[\textbf{output}:] folding automata on $F_{g,n}$  
\end{itemize}
\begin{itemize}
\item[\texttt{step 1}.] Generate the finite set of diffeomorphism types 
of embedded train tracks $\tau_i \subset F_{g,n}$. 
\item[\texttt{step 2}.] 
For each $\tau_i$, compute all the elementary folding maps 
$f_{ij} \colon \tau_i \to \tau'_{ij}$ from $\tau_i$ if any.  
Compute isomorphisms $h_{ij}\colon \tau'_{ij} \to \tau_k$ 
from the train track $\tau'_{ij}$ to one in the set~$\{\tau_i\}$. 
\item[\texttt{step 3}.] For each $\tau_i$, compute the isomorphisms 
$g_{i \ell} \colon \tau_i \to \tau_i$ if any.
\item[\texttt{step 4}.] 
The elementary folding maps $h_{ij}\circ f_{ij}$ 
and the isomorphisms $g_{i \ell}$ form the arrows of the folding 
automata. Compute their transition matrices after labeling 
each edge of all the train tracks $\tau_i$.
\end{itemize}
\end{minipage}
\end{boxedminipage}
\end{center}
Note that for \texttt{step 1} one needs to solve the isomorphism 
problem for embedded train tracks. 
Once \texttt{step 1} is done, implementing the other steps is 
more straightforward.

By running \texttt{genauto}, we obtain 
the folding automata as a collection of connected 
directed graph with each arrow labeled by a transition matrix. 
The goal is to enumerate in the folding automata 
all the closed paths representing pseudo-Anosov mapping classes
with an upper bound for the dilatation.

In this paper we deal with 5-braids using simplified folding automata. 
We ignore infinitesimal edges when computing transition matrices. 
Therefore a closed path in a folding automaton represents 
a pseudo-Anosov braid if and only if  its associated transition matrix 
is Perron-Frobenius.

The following is a pseudo-code for our program \texttt{fbrmin}. 
See~\cite{fbrmin} for details. 
\begin{center}
\begin{boxedminipage}{0.8\textwidth}
\vspace*{0.5\topskip}
\begin{center}\texttt{fbrmin}\end{center}
\hskip1.5em
\begin{minipage}{0.85\textwidth}
\begin{itemize}
\item[\textbf{input}:] a directed graph $\mathcal{G}$ 
with arrows labeled by transition matrices, \\
an upper bound $\lambda$ for the minimum dilatation, \\
and a set $\mathcal{W}$ of 
sub-words which are to be avoided during the search 
\item[\textbf{output}:] the list of closed paths in $\mathcal{G}$ 
representing pseudo-Anosov braids with dilatation less than $\lambda$
\end{itemize}
\begin{itemize}
\item[\texttt{step 1}.] Set $\mathtt{maxnorm} = 
\lfloor \lambda^n + n - 1 \rfloor $ 
where $n$ is the dimension of the transition matrices, \\
and set $\mathtt{archive} = \{ \}$, in which 
closed paths with small dilatation are to be stored. 
\item[\texttt{step 2}.] 
Set $\mathtt{tmp}_1 $ to be the set of length-one paths in $\mathcal{G}$.
\item[\texttt{step 3}.] 
For each $i$ from 2 to $\mathtt{maxnorm}$, 
\begin{enumerate}
\item compute $\mathtt{childrenpaths}_i$ 
by appending paths in $\mathtt{tmp}_1$ 
to paths in $\mathtt{tmp}_{i-1}$, 
in all the ways possible in $\mathcal{G}$
\item compute $\mathtt{tmp}_i$  the subset of $\mathtt{childrenpaths}_i$, 
consisting of paths $\beta$
without any subword from the avoided-word set $\mathcal{W}$, 
with transition matrix $M_\beta$ such that 
$|M_\beta| \le \mathtt{maxnorm}$,  and $M_\beta$ 
has at least one row and one column whose row (column) sum is less than 3.
\item Take the subset $\mathtt{selectedcans}_i$ 
of $\mathtt{tmp}_i$ consisting of 
\emph{closed} paths representing pseudo-Anosov braids with dilatation 
less than $\lambda$, and append it to $\mathtt{archive}$.
\end{enumerate}
\item[\texttt{step 4}.] Return $\mathtt{archive}$ 
$( = \cup_i \mathtt{selectedcans}_i   )$.
\end{itemize}
\end{minipage}
\end{boxedminipage}
\end{center}
In \texttt{step 3}~(2), we use Lemma~\ref{lemma:bound},~\ref{lem:avoidword} 
to trim out much of unnecessary computation (see Remark~\ref{rem:trim}). 

When the row sums of a transition matrix $M_\beta$ all exceed 
$3$, then the spectral radius of~$M_\beta$ is greater than $3$. 
In this case the same holds for every transition matrix of the form 
$M_{\beta \cdot \gamma} = M_\beta M_\gamma$ 
since $M_\gamma \ge P $ for some permutation matrix $P$. 
Therefore as we are looking for transition matrices with 
spectral radius less than 3, 
we can safely remove such paths $\beta\cdot \gamma$ from our consideration, 
as done in \texttt{step 3}~(2).

For computational aspects, the proof of Theorem~\ref{thm:twotrigon} using 
the automaton in Figure~\ref{fig:twotrigon} is the main part which 
consumes most time and memory.  
On a 2.40 GHz machine, it took 1000 seconds of time and 150 
mega-bytes of memory. 
During the search roughly 85000 many matrices were actually tested 
for its largest eigenvalue.

We do not know how far the same kind of computation 
would work for more complicated surfaces. 
We expect that at least 
the case for 6-braids, hence for genus-2 closed surface, 
can be done on a personal computer without too much difficulty.

\bibliographystyle{abbrv}
\bibliography{fbr}

\end{document}